\documentclass[10pt, reqno]{amsart}
\usepackage{amssymb,latexsym}
\usepackage{eucal}
\usepackage{tikz}
\usepackage{comment}

\newtheorem*{Thm}{Theorem}
\newtheorem*{Prb}{Problem}

\newcommand{\e}{\varepsilon}

\newcommand{\A}{\mathcal{A}}
\newcommand{\ph}{\varphi}

\newcommand{\rr}{\mathrm{r}}

\begin{document}

\title[On a conjecture of Imrich and M\"uller]
 {On a conjecture of Imrich and M\"uller}

\author{S.\ V.\ Ivanov}
\address{Department of Mathematics\\
University of Illinois \\
Urbana\\   IL 61801\\ U.S.A.}
 \email{ivanov@illinois.edu}
\thanks{Supported in part by the NSF under grant  DMS 09-01782.}
\subjclass[2010]{Primary 20E05, 20E07, 20F65, 57M07.}

\begin{abstract}
A conjecture of Imrich and M\"uller  on rank of the intersection of subgroups of free groups is disproved.
\end{abstract}
\maketitle

Let $F$ be a free group of finite rank, let ${\rm r}(F)$ denote the rank of $F$ and let $\bar{\rm r}(F) := \max
({\rm r}(F)-1,0)$ be the reduced rank of $F$.
Recall that the Hanna Neumann conjecture  \cite{N1} on subgroups of a free group  claims that  if $H_1$, $H_2$
are finitely generated subgroups of $F$, then $\bar { \rm{r} } (H_1 \cap H_2)
\le \bar{ {\rm r}} (H_1) \bar { \rm r}  (H_2)$.  It was shown by Hanna
Neumann  \cite{N1}  that  $\bar {\rm{r} } (H_1 \cap H_2) \le 2 \bar{ \rm{r}} (H_1)
\bar { \rm{r}} (H_2)$. For relevant discussions, results,  proofs, and generalizations of this conjecture, the reader is referred to
\cite{D}, \cite{D2}, \cite{DIv}, \cite{DIv2}, \cite{Fr}, \cite{I99}, \cite{I01}, \cite{I08}, \cite{I10},  \cite{KM}, \cite{Min1},  \cite{N2}, \cite{St}.

In 1994, Imrich and M\"uller \cite{ImM} proved that if $H_1$, $H_2$ are finitely generated subgroups of a free group $F$ and one of $H_1$, $H_2$ has finite index in the subgroup  $\langle H_1,  H_2 \rangle$ generated by $H_1$, $H_2$, then
\begin{equation}\label{IM}
\bar {\rm{r} } ( \langle H_1, H_2 \rangle) \bar { \rm{r} } (H_1 \cap H_2)
\le \bar{ \rm{r}} (H_1)\bar{ \rm{r}} (H_2)  .
\end{equation}

Imrich and M\"uller also conjectured in \cite{ImM}  that the inequality \eqref{IM}  holds whenever $H_1$, $H_2$
satisfy the following two conditions (C1)--(C2).
\begin{enumerate}
\item[(C1)]   If $K_1 \cap K_2 = H_1 \cap H_2$, where $K_i$ is a free factor of $H_i$,   then $K_i = H_i$, $i=1,2$.
\item[(C2)]  The intersection $H_1 \cap H_2$ contains no nontrivial free factor of  $\langle H_1,  H_2 \rangle$.
    \end{enumerate}

Moreover, Imrich and M\"uller \cite{ImM} remarked that  the condition (C1) alone is not sufficient to imply \eqref{IM} and they discussed an example for which (C1) holds while \eqref{IM} and (C2) fail. The existence of a nontrivial free factor in the intersection
$H_1 \cap H_2$ was instrumental in construction of the example of \cite{ImM}.

Note that the inequality \eqref{IM} would provide a much tighter than the Hanna Neumann conjecture's bound for $\bar { \rm{r} } (H_1 \cap H_2)$ by incorporating the coefficient $\frac 1{\bar {\rm{r} } ( \langle H_1, H_2 \rangle)}$ when conditions (C1)--(C2) are satisfied and this looks quite remarkable. Curiously, the inequality \eqref{IM}  holds true for subspaces $H_1, H_2$  of a finite dimensional vector space $V$ if $ \bar \rr(H)$ is understood as the dimension of a subspace $H \subseteq V$.  Nevertheless, this  conjecture of Imrich and M\"uller  can be disproved as follows.

 \begin{Thm} There are  two subgroups $H_1$, $H_2$ of a  free group $F$ such that  $H_1$, $H_2$ satisfy
 conditions (C1)--(C2)  and $\bar{ \rm{r}}(H_1) =3$,  $\bar{ \rm{r}}(H_2) =5$,
 $\bar{ \rm{r} } (H_1 \cap H_2) = 8$, and $\bar {\rm{r} } ( \langle H_1, H_2 \rangle) =2$. In particular, the inequality \eqref{IM}  fails for $H_1, H_2$.
 \end{Thm}

 {\em Proof of  Theorem.}  Let $\A$ be a set, called an alphabet, and let $ \A^{\pm 1} := \A \cup \A^{-1}$, where $\A^{-1}$ is the set of formal inverses of elements of $\A$, $\A^{-1}$ is disjoint from $\A$. Let $F(\A)$ denote the free group whose free generators are elements of  $\A$.
 A graph $\Gamma$ is called {\em labeled} if $\Gamma$ is equipped
 with a function $\ph : E\Gamma \to \A^{\pm 1}$, where $E\Gamma$ is the set of oriented edges of $\Gamma$, such that $\ph(e^{-1}) = \ph(e)^{-1}$ for every $e \in E\Gamma$. If $p = e_1 \ldots e_\ell$ is a path in  a labeled graph $\Gamma$,
 where $e_1, \ldots, e_\ell \in  E\Gamma$, then  the label $\ph(p)$ of $p$ is the word $\ph(p) := \ph(e_1) \ldots \ph(e_\ell)$ over $\A^{\pm 1}$.   The initial vertex of a path $p$ is denoted $p_-$ and the terminal vertex of $p$ is denoted  $p_+$.

 We now construct some subgroups of a free group by means of their Stallings graphs, for definitions see  \cite{D}, \cite{KM}, \cite{St}.

 Let $k, \ell, m, n$ be some positive integers such that  $m \ge \max(\ell, n)$.  Consider three closed labeled paths $q_a, q_b, q_c$ such that $\ph(q_a) = a^{k \ell}$, $\ph(q_b) =  b^{k m}$,  $\ph(q_c) =  c^{k n}$. A vertex of $q_x$, $x \in \{ a,b,c \}$, at a distance from $(q_x)_-=  (q_x)_+$  that is divisible by $k$, is called a {\em phase} vertex of  $q_x$. We connect some $\ell$ of  phase vertices of  $q_b$ to  the $\ell$   phase vertices of  $q_a$ with edges $f_1, \dots, f_\ell$ so that $(f_i)_- \in q_b$, $(f_i)_+ \in q_a$ are all distinct and $\ph(f_i) = d$, $i= 1, \dots, \ell$.  We also connect some $n$ of  phase vertices of  $q_b$ to  the $n$   phase vertices of  $q_c$ with edges $g_1, \dots, g_n$ so that $(g_i)_- \in q_b$, $(g_i)_+ \in q_c$ are all distinct and $\ph(g_i) = e$, $i= 1, \dots,n$. \ Let $\Gamma(k,\ell,m,n)$
denote a labeled graph  constructed this way, where $\A^{\pm 1} = \{ a^{\pm 1}, b^{\pm 1}, c^{\pm 1}, d^{\pm 1}, e^{\pm 1} \}$, see Fig.~1 where
the graph $\Gamma(1,1,1,1)$,  also denoted  $\Gamma(\bar 1)$, is depicted.
\vskip 2.4mm

\begin{center}
\begin{tikzpicture}[scale=.46]
\draw  (-4.5,4) ellipse (1 and 1);
\draw  (-0.5,4) ellipse (1 and 1);
\draw  (-2.5,1) ellipse (1 and 1);
\draw  (-2.5,2) [fill = black] circle (.07);
\node at (-6.1,4) {$a$};
\node at (1.1,4) {$c$};
\node at (-3.6,2.4) {$d$};
\node at (-1.36,2.4) {$e$};
\draw [-latex](0.5,4.1) --(0.5,3.95);
\draw [-latex](-5.5,3.95) --(-5.5,4.1);
\node at (-2.5,-1.2) {Figure~1};
\node at (-2.5,.6) {$b$};

\draw (-2.5,2) --(-3.8,3.3);
\draw (-2.5,2) --(-1.2,3.3);
\draw [-latex](-2.5,2) --(-3.25,2.75);
\draw [-latex](-2.5,2) --(-1.75,2.75);
\draw [-latex](-2.6,0) --(-2.4,0);
\draw  (-1.2,3.3) [fill = black] circle (.07);
\draw  (-3.8,3.3) [fill = black] circle (.07);
\draw [-latex](-5.5,3.95) --(-5.5,4.1);
\draw [-latex](-5.5,3.95) --(-5.5,4.1);

\node at (2.8,1.5) {$\Gamma(1,1,1,1)= \Gamma(\bar 1)$};
\end{tikzpicture}
\end{center}

\noindent
Note $\Gamma(k,\ell,m,n)$ is not uniquely defined because there are choices for   $(f_j)_-, (g_{j'})_- \in q_b$.
Let $H(\Gamma(k,\ell,m,n))$ denote the subgroup of the free group $F(\A)$ which consists of all the words $\ph(p)$,
where $p$ is a closed path in a fixed graph $\Gamma(k,\ell,m,n)$ such that  $p_-= p_+ = (q_b)_-$. It is immediate from the definitions that
$H(\Gamma(\bar 1)) =  \langle b, dad^{- 1}  , ece^{- 1} \rangle$ and, for every $\Gamma(k,\ell,m,n)$, the subgroup $H(\Gamma(k,\ell,m,n))$  is a subgroup of $H(\Gamma(\bar 1))$. It is easy to check that  $\bar { \rm{r} }(H(\Gamma(k,\ell,m,n)))= \ell +n$.
\smallskip

Consider two such graphs $ \Gamma(7,1,3,2)$,  $ \Gamma(11,2,5,3)$ and associated subgroups
$$
H_1 := H(\Gamma(7,1,3,2)) ,   \quad H_2 := H(\Gamma(11,2,5,3)) .
$$
Denote $\Delta_1 :=  \Gamma(7,1,3,2)$   and  $\Delta_2 :=  \Gamma(11,2,5,3)$.
Since $m_i = \ell_i + n_i$ for the graph $\Delta_i$, $i=1,2$, we may assume that every phase vertex of  $\Delta_1$ and $\Delta_2$  has degree 3. Note that the parameters are chosen so that $\gcd(|q_{x, 1}|,  |q_{x, 2}|) = 1$  for every $x \in \{ a,b,c\}$.  Indeed,  $ |q_{a, 1}| =  7$, $| q_{a, 2}| =22$,  $ | q_{b, 1}| =  21$, $| q_{b, 2}| =55$, $ | q_{c, 1}| =  14$, $| q_{c, 2}| =33$. Now  it is easy to see  that the entire pullback $\Delta_1 \underset{\Gamma(\bar 1)}  \times \Delta_2$ is a graph of the form $ \Gamma(77,2,15,6)$. In particular, this pullback is connected and $\bar { \rm{r} }(H_1 \cap H_2 )= -\chi(\Gamma(77,2,15,6))  =8$. On the other hand, it is immediate  that $\langle H_1, H_2 \rangle = H(\Gamma(\bar 1))$, hence $\bar { \rm{r} }(\langle H_1, H_2 \rangle )=2$.  It remains to show that the subgroups $  H_1, H_2 $ satisfy the conditions (C1)--(C2).

Suppose $K_i$ is a finitely generated subgroup of  $H_i$, $i=1,2$, and
$K_1 \cap K_2 = H_1 \cap H_2$. Let $\Delta(K_i)$ denote the Stallings graph of $K_i$, $i=1,2$. Then there is a canonical locally injective graph map $\Delta(K_i) \to \Gamma(\bar 1)$, $i=1,2$, that factors through $\Delta_i$. Since $K_1 \cap K_2 = H_1 \cap H_2$, it follows that $\Psi := $ $\Delta_1 \underset{\Gamma(\bar 1)}  \times \Delta_2$ is equal to  $\mbox{core}(\Delta(K_1) \underset{\Gamma(\bar 1)}  \times \Delta(K_2) )$ and the projection map $\Psi \to \Delta_i$    factors through     $\Delta(K_i)$, $i=1,2$. Observe that for every vertex $v \in V \Delta_i$ of degree 3, there exists a  vertex $w \in V \Psi$ of degree 3 that maps to $v$. Hence, there is a vertex  $u \in V \Delta(K_i)$ of degree 3 that maps to $v$. This means that $-\chi( \Delta(K_i) )  \ge  -\chi( \Delta_i)$, whence   $\bar { \rm{r} }(K_i) \ge  \bar { \rm{r} }(H_i)$, $i=1,2$.  If, in addition, $K_i$ is a free factor of  $H_i$, we also have  $\bar { \rm{r} }(K_i) \le  \bar { \rm{r} }(H_i)$. Therefore, $\bar { \rm{r} }(K_i) =  \bar { \rm{r} }(H_i)$ and $K_i = H_i$ as desired. The condition (C1) is proven.

Observe that if $p_i$ is a closed reduced path in $\Delta_i$, then the sum of exponents $\e$ on edges $x^\e$ of $p$, where $\ph(x) \in \{ a,b,c\}$, is divisible by $7$ if $i=1$ and is divisible by $11$ if $i=2$.
Since  $\langle H_1, H_2 \rangle =   \langle dad^{-1}, b, ece^{-1} \rangle$, we see that no path $p \in H_i$, where $H_i =  \pi_1( \Delta_i, (q_{b,i})_-)$ is the fundamental group of $\Delta_i$ at $(q_{b,i})_-$, can be a free generator of $H(\Gamma(\bar 1)) = \langle dad^{-1}, b, ece^{-1} \rangle$.  Hence, the condition (C2) is also satisfied for $H_1$, $H_2$, and Theorem is proved. \qed.

\medskip

As above, let $H_1, H_2$ be subgroups of a free group $F$ and let  $S(H_1, H_2)$ denote a set of  representatives of those double cosets $H_1 t H_2$ of $F$, $t \in F$,   that have the property $H_1 \cap t H_2 t^{-1} \ne  \{ 1 \}$. Define
$
\bar \rr(H_1, H_2) := \sum_{s \in S(H_1, H_2)} \bar \rr(H_1\cap s H_2 s^{-1}) .
$

While the inequality \eqref{IM} is false in general, some versions of  \eqref{IM}
might still be true. For instance, a generalization of the inequality \eqref{IM}  with
$\bar {\rm{r} } (H_1, H_2 )$ and  $\bar {\rm{r} } ( \langle H_1, H_2,$ $S(H_1, H_2) \rangle)$
in place of  $\bar { \rm{r} } (H_1 \cap H_2)$ and $\bar {\rm{r} } ( \langle H_1, H_2 \rangle)$, resp.,
holds true whenever one of $H_1, H_2$ has finite index in $\langle H_1, H_2, S(H_1, H_2) \rangle$.
In this more general case, Imrich--M\"uller's arguments \cite{ImM} for the proof of inequality  \eqref{IM}   are essentially retained.
Note that this proof is straightforward and independent of an erroneous lemma of \cite{ImM}, discussed in
detail by Kent \cite{Kent}, see also \cite{I15} for a counterexample and for a repair to this lemma.

\medskip

We conclude by mentioning an interesting special case of the inequality \eqref{IM}.

\begin{Prb} Does the inequality \eqref{IM} hold true under the assumption $\bar { \rm{r} } (H_1 \cap H_2) = \bar{ \rm{r}} (H_1)\bar{ \rm{r}} (H_2) > 0$? Equivalently, does the equality $\bar { \rm{r} } (H_1 \cap H_2) = \bar{ \rm{r}} (H_1)\bar{ \rm{r}} (H_2) > 0$ imply that $\bar {\rm{r} } ( \langle H_1, H_2 \rangle) =1$?
\end{Prb}

We remark that this Problem has a positive solution for available examples.
Curiously,
the equality $\bar { \rm{r} } (H_1, H_2) = \bar{ \rm{r}} (H_1)\bar{ \rm{r}} (H_2) > 0$ need not imply the equality
$\bar {\rm{r} } ( \langle H_1, H_2, S(H_1, H_2)  \rangle) =1$ or even a bound
$\bar {\rm{r} } ( \langle H_1, H_2, S(H_1, H_2)  \rangle) < C$, where $C$ is a constant,
as follows from the following.
\medskip

\noindent
{\bf Example.} Suppose that $H_0$, $H_1$ are proper normal subgroups of finite index of a free group $F(a,b)$ and $F(a,b) = H_0 H_1$.  For instance, one can take normal subgroups  $H_0$, $H_1$ of $F(a,b)$ so that the quotient groups $F(a,b)/H_0$, $F(a,b)/H_1$ have coprime orders.
Then  $ |F(a,b) / H_i| = \bar{ \rm{r}} (H_i)$, $i=0,1$, and
 $ \bar { \rm{r} } (H_0 \cap H_1) = \bar{ \rm{r}} (H_0)  \bar{ \rm{r}} (H_1)$.
Next, let $H_2 := \langle H_0,  c_1 H_0 c_1^{-1}, \dots,  c_{k-1} H_0 c_{k-1}^{-1} \rangle $ be a subgroup of  the free group $F(a,b, c_1, \dots, c_{k-1})$.    Then it is not difficult to check that $\bar { \rm{r} } (H_1, H_2) =  \bar{ \rm{r}} (H_1) \bar{ \rm{r}} (H_2)$
 and   $\langle H_1, H_2, S(H_1, H_2) \rangle = F(a,b, c_1, \dots, c_{k-1})$, which implies that
 $$
 \bar {\rm{r} } ( \langle H_1, H_2, S(H_1, H_2)  \rangle) = k .
 $$
Hence, the reduced rank
$\bar {\rm{r} } ( \langle H_1, H_2$, $S(H_1, H_2)  \rangle)$ could be arbitrarily high
in the situation when  $\bar { \rm{r} } (H_1, H_2) =  \bar{ \rm{r}} (H_1) \bar{ \rm{r}} (H_2) > 0$.

\end{document}